\documentclass[12pt]{article}

\usepackage{amsfonts}
\usepackage{amssymb}

\title{Embeddings of vertex operator algebras associated to
orthogonal affine Lie algebras}
\author{Ozren Per\v{s}e}
\date{}
\begin{document}
\def \Z{\Bbb Z}
\def \C{\Bbb C}
\def \R{\Bbb R}
\def \Q{\Bbb Q}
\def \N{\Bbb N}
\def \tr{{\rm tr}}
\def \span{{\rm span}}
\def \Res{{\rm Res}}
\def \End{{\rm End}}
\def \E{{\rm End}}
\def \Ind {{\rm Ind}}
\def \Irr {{\rm Irr}}
\def \Aut{{\rm Aut}}
\def \Hom{{\rm Hom}}
\def \mod{{\rm mod}}
\def \ann{{\rm Ann}}
\def \<{\langle}
\def \>{\rangle}
\def \t{\tau }
\def \a{\alpha }
\def \e{\epsilon }
\def \l{\lambda }
\def \L{\Lambda }
\def \g{\gamma}
\def \b{\beta }
\def \om{\omega }
\def \o{\omega }
\def \c{\chi}
\def \ch{\chi}
\def \cg{\chi_g}
\def \ag{\alpha_g}
\def \ah{\alpha_h}
\def \ph{\psi_h}
\def \be{\begin{equation}\label}
\def \ee{\end{equation}}
\def \bl{\begin{lem}\label}
\def \el{\end{lem}}
\def \bt{\begin{thm}\label}
\def \et{\end{thm}}
\def \bp{\begin{prop}\label}
\def \ep{\end{prop}}
\def \br{\begin{rem}\label}
\def \er{\end{rem}}
\def \bc{\begin{coro}\label}
\def \ec{\end{coro}}
\def \bd{\begin{de}\label}
\def \ed{\end{de}}
\def \pf{{\bf Proof. }}
\def \voa{{vertex operator algebra}}

\newtheorem{thm}{Theorem}[section]
\newtheorem{prop}[thm]{Proposition}
\newtheorem{coro}[thm]{Corollary}
\newtheorem{conj}[thm]{Conjecture}
\newtheorem{lem}[thm]{Lemma}
\newtheorem{rem}[thm]{Remark}
\newtheorem{de}[thm]{Definition}
\newtheorem{hy}[thm]{Hypothesis}
\makeatletter \@addtoreset{equation}{section}
\def\theequation{\thesection.\arabic{equation}}
\makeatother \makeatletter

\newcommand{\binom}[2]{{{#1}\choose {#2}}}
    \newcommand{\nno}{\nonumber}
    \newcommand{\lbar}{\bigg\vert}
    \newcommand{\p}{\partial}
    \newcommand{\dps}{\displaystyle}
    \newcommand{\bra}{\langle}
    \newcommand{\ket}{\rangle}
 \newcommand{\res}{\mbox{\rm Res}}
\renewcommand{\hom}{\mbox{\rm Hom}}
  \newcommand{\epf}{\hspace{2em}$\Box$}
 \newcommand{\epfv}{\hspace{1em}$\Box$\vspace{1em}}
\newcommand{\nord}{\mbox{\scriptsize ${\circ\atop\circ}$}}
\newcommand{\wt}{\mbox{\rm wt}\ }

\maketitle
\begin{abstract}
Let $L_{D_{\ell}}(-\ell +\frac{3}{2},0)$ (resp. $L_{B_{\ell}}(-\ell
+\frac{3}{2},0)$) be the simple vertex operator algebra associated
to affine Lie algebra of type $D_{\ell}^{(1)}$ (resp.
$B_{\ell}^{(1)}$) with the lowest admissible half-integer level
$-\ell + \frac{3}{2}$. We show that $L_{D_{\ell}}(-\ell
+\frac{3}{2},0)$ is a vertex subalgebra of $L_{B_{\ell}}(-\ell
+\frac{3}{2},0)$ with the same conformal vector. For $\ell =4$,
$L_{D_{4}}(-\frac{5}{2},0)$ is a vertex subalgebra of three copies
of $L_{B_{4}}(-\frac{5}{2},0)$ contained in
$L_{F_{4}}(-\frac{5}{2},0)$, and all five of these vertex operator
algebras have the same conformal vector.
\end{abstract}


\section{Introduction}

Let $\frak g _{X_{l}}$ be the simple Lie algebra of type $X_{l}$,
$\hat{\frak g}_{X_{l}}$ the associated affine Lie algebra of type
$X_{l}^{(1)}$, and $L_{X_{l}}(k,0)$ the simple vertex operator
algebra associated to $\hat{\frak g}_{X_{l}}$ of level $k \in \C$,
$k \neq - h^{\vee}$, with conformal vector obtained from the
Segal-Sugawara construction. It is an interesting problem to find
pairs of vertex operator algebras $L_{X_{l}}(k,0)$ and
$L_{Y_{l'}}(k',0)$, such that $L_{Y_{l'}}(k',0)$ is a vertex
subalgebra of $L_{X_{l}}(k,0)$ with the same conformal vector. In
the case of positive integer levels, such examples are called
conformal embeddings in mathematical physics (cf. \cite{AGO},
\cite{BB}, \cite{SW}).

In \cite{P2}, we showed that $L_{B_{4}}(-\frac{5}{2},0)$ is a vertex
subalgebra of $L_{F_{4}}(-\frac{5}{2},0)$ with the same conformal
vector. The motivation for studying level $k=-\frac{5}{2}$ was that
the corresponding central charge of these vertex operator algebras
is the same (and equal to $-20$). In this case, one can also
consider Lie subalgebra $\frak g _{D_{4}}$ of $\frak g _{B_{4}}$ and
the associated vertex operator algebra $L_{D_{4}}(-\frac{5}{2},0)$.
It turns out that this vertex operator algebra also has central
charge $-20$.

More generally, one can consider the embedding of orthogonal Lie
algebras $\frak g _{D_{l}}$ into $\frak g _{B_{l}}$, for any $l \geq
4$. On the level of affine Lie algebras and corresponding vertex
operator algebras, this embedding has been studied for level $k=1$
(cf. \cite{FF}). Using explicit realizations of modules
$L_{D_{l}}(1,0)$ and $L_{B_{l}}(1,0)$, one can see that
$L_{D_{l}}(1,0)$ is a vertex subalgebra of $L_{B_{l}}(1,0)$, but the
corresponding conformal vectors are different in this case.

In this paper we study a level $k$ such that $L_{D_{l}}(k,0)$ is a
vertex subalgebra of $L_{B_{l}}(k,0)$ with the same conformal
vector. The equality of conformal vectors implies the equality of
corresponding central charges
\begin{eqnarray*}
\frac{k \dim \frak g _{D_{l}}}{k+h_{D_{l}}^{\vee}}= \frac{k \dim
\frak g _{B_{l}}}{k+h_{B_{l}}^{\vee}}
\end{eqnarray*}
Apart from the trivial solution $k=0$, the only solution of this
equation is $k=-l+\frac{3}{2}$. This suggests the study of vertex
operator algebras $L_{D_{l}}(-l+\frac{3}{2},0)$ and
$L_{B_{l}}(-l+\frac{3}{2},0)$. Since there are no known realizations
of these irreducible modules, we consider these modules as quotients
of the corresponding generalized Verma modules.

Vertex operator algebra $L_{B_{l}}(-l+\frac{3}{2},0)$ has been
studied in \cite{P1}. The level $k=-l+\frac{3}{2}$ is an admissible
level for $\hat{\frak g}_{B_{l}}$, in the sense of Kac and Wakimoto
(cf. \cite{KW1}, \cite{KW2}, \cite{W}). The results on admissible
modules from \cite{KW1} and \cite{KW2} imply that
$L_{B_{l}}(-l+\frac{3}{2},0)$ is a quotient of the corresponding
generalized Verma module modulo the ideal generated by one singular
vector $v_{B_{l}}$ of conformal weight $2$. The explicit formula for
that vector was determined in \cite{P1}, and is relatively simple
(cf. relation (\ref{rel.sing.v.B})). Vertex operator algebras
associated to affine Lie algebras with admissible levels have also
recently been studied in \cite{A1,A2,AM,DLM}.

In this paper we study vertex operator algebra
$L_{D_{l}}(-l+\frac{3}{2},0)$. We verify that level
$k=-l+\frac{3}{2}$ is admissible for $\hat{\frak g}_{D_{l}}$, and
show (using \cite{KW1}, \cite{KW2}) that
$L_{D_{l}}(-l+\frac{3}{2},0)$ is a quotient of corresponding
generalized Verma module modulo the ideal generated by one singular
vector $v_{D_{l}}$ of conformal weight $4$. The formula for that
singular vector is much more complicated than the formula for
$v_{B_{l}}$ (cf. relation (\ref{rel.sing.v.D})). Using explicit
formulas for vectors $v_{B_{l}}$ and $v_{D_{l}}$ we show that the
vector $v_{D_{l}}$ is contained in the $\hat{\frak
g}_{B_{l}}$-submodule generated by $v_{B_{l}}$ of the corresponding
generalized Verma module, which directly implies that
$L_{D_{l}}(-l+\frac{3}{2},0)$ is a vertex subalgebra of
$L_{B_{l}}(-l+\frac{3}{2},0)$. Using certain relations generated by
vector $v_{B_{l}}$ on the quotient $L_{B_{l}}(-l+\frac{3}{2},0)$ of
generalized Verma module, we show that the conformal vectors of
$L_{D_{l}}(-l+\frac{3}{2},0)$ and $L_{B_{l}}(-l+\frac{3}{2},0)$ are
the same.

In the special case $l=4$, $L_{F_{4}}(-\frac{5}{2},0)$ contains
three copies of $L_{B_{4}}(-\frac{5}{2},0)$ as vertex subalgebras
with the same conformal vector (cf. \cite{P2}). It follows that
vertex operator algebra $L_{D_{4}}(-\frac{5}{2},0)$ is a vertex
subalgebra of all three copies, with the same conformal vector.

\section{Preliminaries}
\label{sect.prelim}

Let $(V, Y, {\bf 1}, \omega)$ be a vertex operator algebra (cf.
\cite{Bo,FHL,FLM,LL}). A {\it vertex subalgebra} of vertex algebra
$V$ is a subspace $U$ of $V$ such that ${\bf 1} \in U$ and $ Y(a,z)U
\subseteq U[[z,z^{-1}]] $ for any $a \in U$. Suppose that $(V, Y,
{\bf 1}, \omega)$ is a vertex operator algebra and $(U, Y, {\bf 1},
\omega ')$ a vertex subalgebra of $V$, that has a structure of
vertex operator algebra. We say that $U$ is a {\it vertex operator
subalgebra} of $V$ if $\omega ' = \omega$.

Let ${\frak g}$ be a simple Lie algebra over ${\C}$ with a
triangular decomposition  ${\frak g}={\frak n}_{-} \oplus {\frak h}
\oplus {\frak n}_{+}$. Let $(\cdot, \cdot): {\frak g}\times {\frak
g}\to {\C}$ be the Killing form, normalized by the condition
$(\theta, \theta)=2$, where $\theta$ is the highest root of $\frak
g$. The affine Lie algebra $\hat{\frak g}$ associated to ${\frak g}$
is the vector space ${\frak g}\otimes {\C}[t, t^{-1}] \oplus {\C}c $
equipped with the usual bracket operation and the canonical central
element~$c$ (cf. \cite{K1}). Let $h^{\vee}$ be the dual Coxeter
number of $\hat{\frak g}$. Let $\hat{\frak g}=\hat{\frak n}_{-}
\oplus \hat{\frak h} \oplus \hat{\frak n}_{+}$ be the corresponding
triangular decomposition of $\hat{\frak g}$.

For every weight $\lambda\in \hat{{\frak h}}^{*}$, denote by $M(
\lambda)$ the Verma module for $\hat{\frak g}$ with highest weight
$\lambda$, and by $L( \lambda)$ the irreducible $\hat{\frak
g}$-module with highest weight~$\lambda$. Let $U$ be a ${\frak
g}$--module, and let $k\in {\C}$. Denote by $N(k,U)$ the generalized
Verma $\hat{\frak g}$--module of level $k$ induced from $U$. For
$\mu \in {\frak h}^{*}$, denote by $V(\mu)$ the irreducible highest
weight ${\frak g}$--module with highest weight $\mu$. We shall use
the notation $N(k, \mu)$ to denote the $\hat{\frak g}$--module
$N(k,V(\mu))$. Denote by $J(k, \mu)$ the maximal proper submodule of
$N(k, \mu)$ and by $L(k, \mu)=N(k, \mu) /J(k, \mu)$ the
corresponding irreducible $\hat{\frak g}$--module.

Let $\hat{\Delta}^{\vee \mbox{\scriptsize{re}}}$\ (resp.
$\hat{\Delta}^{\vee \mbox{\scriptsize{re}}}_{+}$) $\subset
\hat{\frak h}$ be the set of real (resp. positive real) coroots of $
\hat{\frak g}$. Fix $\lambda \in \hat{\frak h}^*$. Let
$\hat{\Delta}^{\vee \mbox{\scriptsize{re}}}_{\lambda}= \{\alpha\in
\hat{\Delta}^{\vee \mbox{\scriptsize{re}}} \ |\
\langle\lambda,\alpha\rangle\in{\Z} \}$, $\hat{\Delta}^{\vee
\mbox{\scriptsize{re}}}_{\lambda +}=\hat{\Delta}^{\vee
\mbox{\scriptsize{re}}}_{\lambda} \cap \hat{\Delta}^{\vee
\mbox{\scriptsize{re}}}_{+}$, $\hat{\Pi}^{\vee}$ the set of simple
coroots in $\hat{\Delta}^{\vee \mbox{\scriptsize{re}}}$ and
$\hat{\Pi}^{\vee}_{\lambda }= \{\alpha \in \hat{\Delta}^{\vee
\mbox{\scriptsize{re}}}_{\lambda +}\ \vert \ \alpha$ not equal to a
sum of several coroots from $\hat{\Delta}^{\vee
\mbox{\scriptsize{re}}}_{\lambda +} \}$. Recall that a weight
$\lambda \in \hat{\frak h} ^* $ is called {\it admissible} (cf.
\cite{KW1}, \cite{KW2} and \cite{W}) if the following properties are
satisfied:
\begin{eqnarray*}
& &\langle\lambda + \rho,\alpha\rangle \notin -{\Z}_+ \mbox{ for all
}
\alpha \in \hat{\Delta}^{\vee \mbox{\scriptsize{re}}}_{+}, \\
& &{\Q} \hat{\Delta}^{\vee \mbox{\scriptsize{re}}}_{\lambda}={\Q}
\hat{\Pi}^{\vee} .
\end{eqnarray*}

We shall use the following result of V. Kac and M. Wakimoto:
\begin{prop}[{\cite[Corollary 2.1]{KW1}}] \label{t.KW1}
Let $\lambda $ be an admissible weight. Then
\[
L(\lambda) \cong \frac {M(\lambda)}{\sum_{\alpha \in
\hat{\Pi}^{\vee}_{\lambda }}U ( \hat{\frak g}) v^{\alpha}}\ ,
\]
where $v^{\alpha}\in M(\lambda)$ is a singular vector of weight
$r_{\alpha}. \lambda$, the highest weight vector of
$M(r_{\alpha}.\lambda)=\ U(\hat{\frak g}) v^{\alpha}\subset
M(\lambda)$.
\end{prop}

For $k \neq - h^{\vee}$, $N(k,0)$ has the structure of vertex
operator algebra (cf. \cite{FZ,FrB,K2,LL,L}) with conformal vector
\begin{eqnarray} \label{rel.Virasoro}
\omega=\frac{1}{2(k+h^{\vee})}\sum_{i=1}^{\dim {\frak g}}
a^{i}(-1)b^{i}(-1){\bf 1},
\end{eqnarray}
where $\{a^{i}\}_{i=1, \dots, \dim {\frak g}}$ is an arbitrary basis
of ${\frak g}$, and $\{b^{i}\}_{i=1, \dots, \dim {\frak g}}$  the
corresponding dual basis of ${\frak g}$ with respect to the form
$(\cdot, \cdot)$.

Since every $\hat{\frak g}$-submodule of $N(k,0)$ is also an ideal
in the vertex operator algebra $N(k,0)$, it follows that $L(k,0)$ is
a vertex operator algebra, for any $k\ne -h^{\vee}$.

\section{Simple Lie algebras of type $D_{l}$ and $B_{l}$}
\label{sect.simple}

Let $l \geq 4$ and consider two $l$--dimensional vector spaces $A_1
= \oplus _{i=1} ^{l}{\C} a_i$, $ A_2 = \oplus _{i=1}^ {l} {\C}
a_i^*$. Let $A=A_1 \oplus A_2$. The Clifford algebra $Cliff(A)$ is
the associative algebra over ${\C}$ generated by $A$ and relations
\[
[a_i,a_j]_{+}=[a_i^*,a_j^*]_{+}=0 ,\quad [a_i,a_j^*]_{+}=\delta_{ij}
, \qquad i,j = 1,\dots ,l.
\]
The normal ordering on $A$ is defined by
\[
:\!xy\!:  \ = \frac 1 2 (xy-yx),\ \qquad x,y \in A .
\]
Then (cf. \cite{Bou} and \cite{FF}) all normally ordered quadratic
elements $:\!xy\!:$ span a Lie algebra $\frak g _{D_{l}}$ of type
$D_{l}$ with Cartan subalgebra $\frak h$ spanned by
\[
H_i=:\!a_ia_i^*\!:,  \qquad i=1, \ldots ,l.
\]
Let $\{\epsilon_i\ | \ 1\leq i\leq l \}\subset {\frak h}^* $ be the
dual basis such that $\epsilon_i(H_j)=\delta_{ij}$. The root system
of $\frak g _{D_{l}}$ is given by
\[
\Delta _{D_{l}}=\{ \pm(\epsilon_i \pm \epsilon_j) \ \vert \ 1 \leq
i,j \leq l, i < j \}
\]
with
$\alpha_1=\epsilon_1-\epsilon_2,...,\alpha_{l-1}=\epsilon_{l-1}-
\epsilon_{l}, \alpha_{l}=\epsilon_{l-1}+ \epsilon_{l}$ being a set
of simple roots. The highest root is $\theta=\epsilon_{1}+
\epsilon_{2}$. The set of positive roots is $\Delta_{D_{l}}^{+}=\{
\epsilon_{i} \pm \epsilon_{j} \, \vert \ i<j \}$. We fix the root
vectors
\[
e_{\epsilon_i-\epsilon_j}=\ :\!a_ia_j^*\!:, \quad
e_{\epsilon_i+\epsilon_j}=\ :\!a_ia_j\!:, \quad
f_{\epsilon_i-\epsilon_j}=\ :\!a_ja_i^*\!:, \quad
f_{\epsilon_i+\epsilon_j}=\ :\!a_j^*a_i^*\!:,
\]
for $i,j = 1, \dots ,l$ and $i<j$.

If the linear space $A$ is added on to quadratic elements
$:\!xy\!:$, $x,y \in A$, one obtains the simple Lie algebra $\frak g
_{B_{l}}$ of type $B_{l}$. The root system of $\frak g _{B_{l}}$
with respect to the same Cartan subalgebra $\frak h$ is
\[
\Delta _{B_{l}}= \Delta _{D_{l}} \cup \{ \pm \epsilon_i \ \vert \ 1
\leq i \leq l \}
\]
with $\beta_1=\epsilon_1-\epsilon_2,...,\beta_{l-1}=\epsilon_{l-1}-
\epsilon_{l}, \beta_{l}= \epsilon_{l}$ being a set of simple roots.
The set of positive roots is $\Delta_{B_{l}}^{+}= \Delta_{D_{l}}^{+}
\cup \{ \epsilon_{i} \, \vert \ 1 \leq i \leq l \}$. Fix the root
vectors for short roots
\[
e_{\epsilon_i}= a_i, \quad f_{\epsilon_i}=a_i^*, \quad i = 1, \dots
,l.
\]

Denote by $h_{\alpha}= \alpha ^{\vee}= [e_{\alpha},f_{\alpha}]$
coroots, for any positive root $\alpha \in \Delta_{B_{l}}^{+}$.
Clearly $h_{\epsilon_i}=2 H_i$, for $i=1, \ldots ,l$.

\section{Vertex operator algebras $L_{B_{l}}(-l+\frac{3}{2},0)$ and
$L_{D_{l}}(-l+\frac{3}{2},0)$}

Let $\hat{\frak g}_{D_{l}}$ and $\hat{\frak g}_{B_{l}}$ be affine
Lie algebras associated to $\frak g _{D_{l}}$ and  $\frak g
_{B_{l}}$, respectively. For $k \in \C$, denote by $N_{D_{l}}(k,0)$
and $N_{B_{l}}(k,0)$ generalized Verma modules associated to
$\hat{\frak g}_{D_{l}}$ and $\hat{\frak g}_{B_{l}}$ of level $k$,
and by $L_{D_{l}}(k,0)$ and $L_{B_{l}}(k,0)$ corresponding
irreducible modules.

We shall also use the notation $M_{D_{l}}( \lambda)$ (resp.
$M_{B_{l}}( \lambda)$) for the Verma module for $\hat{\frak
g}_{D_{l}}$ (resp. $\hat{\frak g}_{B_{l}}$) with highest weight
$\lambda\in \hat{{\frak h}}^{*}$, and $L_{D_{l}}( \lambda)$ (resp.
$L_{B_{l}}( \lambda)$) for the irreducible $\hat{\frak
g}_{D_{l}}$-module (resp. $\hat{\frak g}_{B_{l}}$-module) with
highest weight $\lambda\in \hat{{\frak h}}^{*}$.

Since $\frak g _{D_{l}}$ is a Lie subalgebra of $\frak g _{B_{l}}$,
it follows that $\hat{\frak g}_{D_{l}}$ is a Lie subalgebra of
$\hat{\frak g}_{B_{l}}$. Using the P-B-W theorem (which gives an
embedding of the universal enveloping algebra of a Lie subalgebra
into the universal enveloping algebra of a Lie algebra) we obtain an
embedding of generalized Verma module $N_{D_{l}}(k,0)$ into
$N_{B_{l}}(k,0)$. Therefore, $N_{D_{l}}(k,0)$ is a vertex subalgebra
of $N_{B_{l}}(k,0)$.

\subsection{Vertex operator algebra $L_{B_{l}}(-l+\frac{3}{2},0)$}

In this subsection we recall some facts about vertex operator
algebra $L_{B_{l}}(-l+\frac{3}{2},0)$. It was proved in
\cite[Theorem 11]{P1} that the maximal $\hat{\frak
g}_{B_{l}}$-submodule of $N_{B_{l}}(-l+\frac{3}{2},0)$ is generated
by one singular vector:

\begin{prop} \label{p.3.1} The maximal $\hat{\frak g}_{B_{l}}$-submodule
of $N_{B_{l}}(-l+\frac{3}{2},0)$ is $J_{B_{l}}(-l+\frac{3}{2},0)=
U(\hat{\frak g}_{B_{l}})v_{B_{l}}$, where
\begin{eqnarray} \label{rel.sing.v.B}
&& v_{B_{l}}= -\frac{1}{4}e_{\epsilon_1}(-1)^2 {\bf 1}+
\sum_{j=2}^{l} e_{\epsilon_1 - \epsilon_j}(-1) e_{\epsilon_1 +
\epsilon_j}(-1){\bf 1}
\end{eqnarray}
is a singular vector in $N_{B_{l}}(-l+\frac{3}{2},0)$.
\end{prop}

Thus
\[
L_{B_{l}}(-l+\frac{3}{2},0) \cong \frac{N_{B_{l}}(-l+\frac{3}{2},0)}
{U(\hat{\frak g}_{B_{l}})v_{B_{l}}}.
\]

\subsection{Vertex operator algebra $L_{D_{l}}(-l+\frac{3}{2},0)$}

In this subsection we show that the maximal $\hat{\frak
g}_{D_{l}}$-submodule of $N_{D_{l}}(-l+\frac{3}{2},0)$ is generated
by one singular vector. We need two lemmas to prove that.

Denote by $\lambda$ the weight $(-l+\frac{3}{2})\Lambda_{0}$. Then
$N_{D_{l}}(-l+\frac{3}{2},0)$ is a quotient of the Verma module
$M_{D_{l}}(\lambda)$ and $L_{D_{l}}(-l+\frac{3}{2},0) \cong
L_{D_{l}}(\lambda)$.

\begin{lem} \label{l.3.4}
The weight $\lambda =(-l+\frac{3}{2})\Lambda_{0}$ is admissible for
$\hat{\frak g}_{D_{l}}$ and
\[
\hat{\Pi}^{\vee}_{\lambda}= \{ (2 \delta - \theta)^{\vee},
\alpha_{1}^{\vee}, \ldots, \alpha_{l}^{\vee} \}.
\]
\end{lem}
{\bf Proof:}
Clearly
\begin{eqnarray*}
&& \langle \lambda + \rho,\alpha _{i}^{\vee}\rangle = 1 \ \ \mbox{for } 1 \leq i \leq l, \\
&& \langle \lambda + \rho,\alpha _{0}^{\vee} \rangle
=-l+\frac{5}{2},
\end{eqnarray*}
which implies
\begin{eqnarray*}
\langle \lambda + \rho,(2 \delta - \theta)^{\vee} \rangle = \langle
\lambda + \rho, 2 \alpha _{0}^{\vee}+ \alpha _{1}^{\vee}+ 2 \alpha
_{2}^{\vee}+ \ldots 2 \alpha _{l-2}^{\vee} +  \alpha _{l-1}^{\vee}+
\alpha _{l}^{\vee} \rangle = 2.
\end{eqnarray*}
The claim of lemma now follows easily.
$\;\;\;\;\Box$

\begin{lem} \label{l.3.5}
Vector
\begin{eqnarray}
&& v_{D_{l}}= \frac{2(2l+1)}{3} \sum_{i=3}^{l} \sum_{j=3 \atop j
\neq i}^{l} e_{\epsilon_1 - \epsilon_i}(-1) e_{\epsilon_1 +
\epsilon_i}(-1) e_{\epsilon_2 - \epsilon_j}(-1) e_{\epsilon_2 +
\epsilon_j}(-1) {\bf 1} \nonumber \\
&& + \frac{2l+1}{3} \sum_{i=3}^{l} e_{\epsilon_1 - \epsilon_i}(-1)
e_{\epsilon_1 + \epsilon_i}(-1) e_{\epsilon_2 - \epsilon_i}(-1)
e_{\epsilon_2 +
\epsilon_i}(-1) {\bf 1} \nonumber \\
&& -\frac{2l+1}{3} \sum_{i=3}^{l} \sum_{j=3 \atop j \neq i}^{l}
e_{\epsilon_1 - \epsilon_i}(-1) e_{\epsilon_2 + \epsilon_i}(-1)
e_{\epsilon_1 + \epsilon_j}(-1) e_{\epsilon_2 -
\epsilon_j}(-1) {\bf 1} \nonumber \\
&& -\frac{2l+1}{6} \sum_{i=3}^{l} \sum_{j=3}^{l} e_{\epsilon_1 -
\epsilon_i}(-1) e_{\epsilon_2 + \epsilon_i}(-1) e_{\epsilon_1 -
\epsilon_j}(-1) e_{\epsilon_2 +
\epsilon_j}(-1) {\bf 1} \nonumber \\
&& -\frac{2l+1}{6} \sum_{i=3}^{l} \sum_{j=3}^{l} e_{\epsilon_1 +
\epsilon_i}(-1) e_{\epsilon_2 - \epsilon_i}(-1) e_{\epsilon_1 +
\epsilon_j}(-1) e_{\epsilon_2 -
\epsilon_j}(-1) {\bf 1} \nonumber \\
&& +2 \sum_{i=3}^{l} \sum_{j=3}^{i-1} e_{\epsilon_1 +
\epsilon_2}(-1) e_{\epsilon_1 + \epsilon_j}(-1) e_{\epsilon_2 -
\epsilon_i}(-1) f_{\epsilon_j - \epsilon_i}(-1) {\bf 1} \nonumber \\
&& +2 \sum_{i=3}^{l} \sum_{j=i+1}^{l} e_{\epsilon_1 +
\epsilon_2}(-1) e_{\epsilon_1 + \epsilon_j}(-1) e_{\epsilon_2 -
\epsilon_i}(-1) e_{\epsilon_i - \epsilon_j}(-1) {\bf 1} \nonumber \\
&& -2 \sum_{i=3}^{l} \sum_{j=3}^{i-1} e_{\epsilon_1 +
\epsilon_2}(-1) e_{\epsilon_1 - \epsilon_j}(-1) e_{\epsilon_2 -
\epsilon_i}(-1) e_{\epsilon_j + \epsilon_i}(-1) {\bf 1} \nonumber \\
&& +2 \sum_{i=3}^{l} \sum_{j=i+1}^{l} e_{\epsilon_1 +
\epsilon_2}(-1) e_{\epsilon_1 - \epsilon_j}(-1) e_{\epsilon_2 -
\epsilon_i}(-1) e_{\epsilon_i + \epsilon_j}(-1) {\bf 1} \nonumber \\
&& +2 \sum_{i=3}^{l} \sum_{j=3}^{i-1} e_{\epsilon_1 +
\epsilon_2}(-1) e_{\epsilon_2 + \epsilon_i}(-1) e_{\epsilon_1 +
\epsilon_j}(-1) f_{\epsilon_j + \epsilon_i}(-1) {\bf 1} \nonumber \\
&& -2 \sum_{i=3}^{l} \sum_{j=i+1}^{l} e_{\epsilon_1 +
\epsilon_2}(-1) e_{\epsilon_2 + \epsilon_i}(-1) e_{\epsilon_1 +
\epsilon_j}(-1) f_{\epsilon_i + \epsilon_j}(-1) {\bf 1} \nonumber \\
&& -2 \sum_{i=3}^{l} \sum_{j=3}^{i-1} e_{\epsilon_1 +
\epsilon_2}(-1) e_{\epsilon_2 + \epsilon_i}(-1) e_{\epsilon_1 -
\epsilon_j}(-1) e_{\epsilon_j - \epsilon_i}(-1) {\bf 1} \nonumber \\
&& -2 \sum_{i=3}^{l} \sum_{j=i+1}^{l} e_{\epsilon_1 +
\epsilon_2}(-1) e_{\epsilon_2 + \epsilon_i}(-1) e_{\epsilon_1 -
\epsilon_j}(-1) f_{\epsilon_i - \epsilon_j}(-1) {\bf 1} \nonumber \\
&& - \frac{2(2l-5)}{3} \sum_{i=3}^{l} e_{\epsilon_1 +
\epsilon_2}(-1) f_{\epsilon_1 - \epsilon_2}(-1) e_{\epsilon_1 -
\epsilon_i}(-1) e_{\epsilon_1 +
\epsilon_i}(-1) {\bf 1} \nonumber \\
&& + \frac{2(2l-5)}{3} \sum_{i=3}^{l} e_{\epsilon_1 +
\epsilon_2}(-1) e_{\epsilon_1 - \epsilon_2}(-1) e_{\epsilon_2 -
\epsilon_i}(-1) e_{\epsilon_2 +
\epsilon_i}(-1) {\bf 1} \nonumber \\
&& + \sum_{i=3}^{l} e_{\epsilon_1 + \epsilon_2}(-1) e_{\epsilon_1 +
\epsilon_i}(-1) e_{\epsilon_2 - \epsilon_i}(-1) h_{\epsilon_i}(-1) {\bf 1} \nonumber \\
&& + \frac{2l-5}{3} \sum_{i=3}^{l} e_{\epsilon_1 + \epsilon_2}(-1)
e_{\epsilon_1 +
\epsilon_i}(-1) e_{\epsilon_2 - \epsilon_i}(-1) h_{\epsilon_1 - \epsilon_2}(-1) {\bf 1} \nonumber \\
&& - \sum_{i=3}^{l} e_{\epsilon_1 + \epsilon_2}(-1) e_{\epsilon_1 -
\epsilon_i}(-1) e_{\epsilon_2 + \epsilon_i}(-1) h_{\epsilon_i}(-1) {\bf 1} \nonumber \\
&& + \frac{2l-5}{3} \sum_{i=3}^{l} e_{\epsilon_1 + \epsilon_2}(-1)
e_{\epsilon_1 -
\epsilon_i}(-1) e_{\epsilon_2 + \epsilon_i}(-1) h_{\epsilon_1 - \epsilon_2}(-1) {\bf 1} \nonumber \\
&& + \frac{2l-5}{3} \sum_{i=3}^{l} e_{\epsilon_1 + \epsilon_2}(-1)
e_{\epsilon_1 +
\epsilon_i}(-2) e_{\epsilon_2 - \epsilon_i}(-1) {\bf 1} \nonumber \\
&& + \frac{2l-5}{3} \sum_{i=3}^{l} e_{\epsilon_1 + \epsilon_2}(-1)
e_{\epsilon_1 -
\epsilon_i}(-2) e_{\epsilon_2 + \epsilon_i}(-1) {\bf 1} \nonumber \\
&& - \frac{(2l+1)(l-3)}{3} \sum_{i=3}^{l} e_{\epsilon_1 +
\epsilon_2}(-2) e_{\epsilon_1 +
\epsilon_i}(-1) e_{\epsilon_2 - \epsilon_i}(-1) {\bf 1} \nonumber \\
&& - \frac{(2l+1)(l-3)}{3} \sum_{i=3}^{l} e_{\epsilon_1 +
\epsilon_2}(-2) e_{\epsilon_1 -
\epsilon_i}(-1) e_{\epsilon_2 + \epsilon_i}(-1) {\bf 1} \nonumber \\
&& - \frac{2l-5}{3} \sum_{i=3}^{l} e_{\epsilon_1 + \epsilon_2}(-1)
e_{\epsilon_1 +
\epsilon_i}(-1) e_{\epsilon_2 - \epsilon_i}(-2) {\bf 1} \nonumber \\
&& - \frac{2l-5}{3} \sum_{i=3}^{l} e_{\epsilon_1 + \epsilon_2}(-1)
e_{\epsilon_1 -
\epsilon_i}(-1) e_{\epsilon_2 + \epsilon_i}(-2) {\bf 1} \nonumber \\
&& + \frac{(2l-5)(2l-1)}{6}  e_{\epsilon_1 + \epsilon_2}(-2)
e_{\epsilon_1 +
\epsilon_2}(-1) h_{\epsilon_1 - \epsilon_2}(-1) {\bf 1} \nonumber \\
&& - \frac{2l-5}{2} e_{\epsilon_1 + \epsilon_2}(-2) e_{\epsilon_1
+ \epsilon_2}(-1) h_{\epsilon_1}(-1) {\bf 1} \nonumber \\
&& - \frac{(2l+1)(2l-5)(2l-7)}{24} e_{\epsilon_1 + \epsilon_2}(-2)^2
{\bf 1} + \frac{(2l-5)^2}{2} e_{\epsilon_1 + \epsilon_2}(-1)
e_{\epsilon_1 +
\epsilon_2}(-3) {\bf 1} \nonumber \\
&& - \frac{2(2l-5)(l-2)}{3(2l-1)}e_{\epsilon_1 + \epsilon_2}(-1)^2
(e_{\epsilon_1 - \epsilon_2}(-1)f_{\epsilon_1 - \epsilon_2}(-1)+
f_{\epsilon_1 - \epsilon_2}(-1)e_{\epsilon_1 - \epsilon_2}(-1)){\bf
1} \nonumber\\
&& + \frac{2l-5}{2l-1} \sum_{i=3}^{l}e_{\epsilon_1 +
\epsilon_2}(-1)^2 (e_{\epsilon_2 - \epsilon_i}(-1)f_{\epsilon_2 -
\epsilon_i}(-1)+ f_{\epsilon_2 - \epsilon_i}(-1)e_{\epsilon_2 -
\epsilon_i}(-1)){\bf 1} \nonumber\\
&& + \frac{2l-5}{2l-1} \sum_{i=3}^{l}e_{\epsilon_1 +
\epsilon_2}(-1)^2 (e_{\epsilon_2 + \epsilon_i}(-1)f_{\epsilon_2 +
\epsilon_i}(-1)+ f_{\epsilon_2 + \epsilon_i}(-1)e_{\epsilon_2 +
\epsilon_i}(-1)){\bf 1} \nonumber\\
&& + \frac{2l-5}{2l-1} \sum_{i=3}^{l}e_{\epsilon_1 +
\epsilon_2}(-1)^2 (e_{\epsilon_1 - \epsilon_i}(-1)f_{\epsilon_1 -
\epsilon_i}(-1)+ f_{\epsilon_1 - \epsilon_i}(-1)e_{\epsilon_1 -
\epsilon_i}(-1)){\bf 1} \nonumber\\
&& + \frac{2l-5}{2l-1} \sum_{i=3}^{l}e_{\epsilon_1 +
\epsilon_2}(-1)^2 (e_{\epsilon_1 + \epsilon_i}(-1)f_{\epsilon_1 +
\epsilon_i}(-1)+ f_{\epsilon_1 + \epsilon_i}(-1)e_{\epsilon_1 +
\epsilon_i}(-1)){\bf 1} \nonumber\\
&& + \frac{2l-5}{2l-1} e_{\epsilon_1 + \epsilon_2}(-1)^2
(e_{\epsilon_1 + \epsilon_2}(-1)f_{\epsilon_1 + \epsilon_2}(-1)+
f_{\epsilon_1 + \epsilon_2}(-1)e_{\epsilon_1 +
\epsilon_2}(-1)){\bf 1} \nonumber\\
&& - \frac{4}{2l-1} \sum_{\alpha \in \Delta _{D_{l}}^{+} \atop (
\alpha , \epsilon_1)=0 , ( \alpha , \epsilon_2)=0 }e_{\epsilon_1 +
\epsilon_2}(-1)^2( e_{\alpha}(-1) f_{\alpha}(-1){\bf 1}
+f_{\alpha}(-1)e_{\alpha}(-1) ){\bf 1} \nonumber \\
&& +e_{\epsilon_1 + \epsilon_2}(-1)^2 \Big( -\frac{2l-5}{6}
h_{\epsilon_1 - \epsilon_2}(-1)^2 {\bf 1} - \frac{1}{2(2l-1)}
h_{\epsilon_1}(-1)^2 {\bf 1} \nonumber \\
&& - \frac{2l-5}{2l-1}h_{\epsilon_1 - \epsilon_2}(-1)h_{\epsilon_1 +
\epsilon_2}(-1){\bf 1} + \frac{4}{2l-1} \sum_{i=3}^{l} h_{\epsilon_1
- \epsilon_i}(-1)h_{\epsilon_1 + \epsilon_i}(-1){\bf 1} \Big) \nonumber \\
&& + \frac{2l-5}{2} e_{\epsilon_1 + \epsilon_2}(-1)^2 h_{\epsilon_1
+ \epsilon_2}(-2){\bf 1} \label{rel.sing.v.D}
\end{eqnarray} is a singular vector in
$N_{D_{l}}(-l+\frac{3}{2},0)$.
\end{lem}
{\bf Proof:} Direct verification of relations $e_{\epsilon_k -
\epsilon_{k+1}}(0).v_{D_{l}}=0$, for $k=1,\ldots ,l-1$,
$e_{\epsilon_{l-1} + \epsilon_{l}}(0).v_{D_{l}}=0$ and
$f_{\theta}(1).v_{D_{l}}=0$ (see Appendix).  $\;\;\;\;\Box$

\begin{thm} \label{t.3.6} The maximal $\hat{\frak g}_{D_{l}}$-submodule
of $N_{D_{l}}(-l+\frac{3}{2},0)$ is $J_{D_{l}}(-l+\frac{3}{2},0)=
U(\hat{\frak g}_{D_{l}})v_{D_{l}}$, where $v_{D_{l}}$ is given by
relation (\ref{rel.sing.v.D}).
\end{thm}
{\bf Proof:} It follows from Proposition \ref{t.KW1} and Lemma
\ref{l.3.4} that the maximal submodule of the Verma module
$M_{D_{l}}(\lambda)$ is generated by $l+1$ singular vectors with
weights
\[
r_{2 \delta - \theta}.\lambda, \ r_{\alpha_{1}}.\lambda, \ldots,
r_{\alpha_{l}}.\lambda.
\]
It follows from Lemma \ref{l.3.5} that $v_{D_{l}}$ is a singular
vector with weight $\lambda -4\delta +2 \theta=r_{2\delta -
\theta}.\lambda$. Other singular vectors have weights
\[
r_{\alpha_{i}}.\lambda =\lambda - \langle \lambda + \rho,
\alpha_{i}^{\vee} \rangle \alpha_{i}=\lambda - \alpha_{i}, \ 1 \leq
i \leq l,
\]
so the images of these vectors under the projection of
$M_{D_{l}}(\lambda)$ onto $N_{D_{l}}(-l+\frac{3}{2},0)$ are 0.
Therefore, the maximal submodule of $N_{D_{l}}(-l+\frac{3}{2},0)$ is
generated by the vector $v_{D_{l}}$, i.e.
$J_{D_{l}}(-l+\frac{3}{2},0)= U(\hat{\frak g}_{D_{l}})v_{D_{l}}$.
$\;\;\;\;\Box$

It follows that
\[
L_{D_{l}}(-l+\frac{3}{2},0) \cong \frac{N_{D_{l}}(-l+\frac{3}{2},0)}
{U(\hat{\frak g}_{D_{l}})v_{D_{l}}}.
\]

\section{Embedding of $L_{D_{l}}(-l+\frac{3}{2},0)$ into  $L_{B_{l}}(-l+\frac{3}{2},0)$}

In this section we show that vector $v_{D_{l}}$ is in the
$\hat{\frak g}_{B_{l}}$-submodule $J_{B_{l}}(-l+\frac{3}{2},0)$ of
$N_{B_{l}}(-l+\frac{3}{2},0)$ generated by vector $v_{B_{l}}$:

\begin{lem} \label{l.emb.tech}
$v_{D_{l}} \in J_{B_{l}}(-l+\frac{3}{2},0)$
\end{lem}
{\bf Proof:} Using relation (\ref{rel.sing.v.B}), one can show that
the following relations hold in $J_{B_{l}}(-l+\frac{3}{2},0)$:
\begin{eqnarray}
&& f_{\epsilon_1}(0). v_{B_{l}}=\frac{1}{2} e_{\epsilon_1}(-1)
h_{\epsilon_1}(-1) {\bf 1} +\sum_{j=2}^{l} e_{\epsilon_1 +
\epsilon_j}(-1) f_{\epsilon_j}(-1) {\bf 1} + \sum_{j=2}^{l}
e_{\epsilon_1 - \epsilon_j}(-1) e_{\epsilon_j}(-1) {\bf 1} \nonumber \\
&& \quad -\frac{2l-3}{2}e_{\epsilon_1}(-2) {\bf 1}, \label{rel.emb.1} \\
&& f_{\epsilon_1 - \epsilon_i}(0). v_{B_{l}}=-\frac{1}{2}
e_{\epsilon_1}(-1) e_{\epsilon_i}(-1) {\bf 1} +\sum_{j=2}^{i-1}
e_{\epsilon_1 + \epsilon_j}(-1) f_{\epsilon_j - \epsilon_i}(-1) {\bf
1} \nonumber \\
&& \quad - \sum_{j=2}^{i-1} e_{\epsilon_1 - \epsilon_j}(-1)
e_{\epsilon_j + \epsilon_i}(-1) {\bf 1} - h_{\epsilon_1 -
\epsilon_i}(-1) e_{\epsilon_1 + \epsilon_i}(-1) {\bf 1}+
\sum_{j=i+1}^{l} e_{\epsilon_1 + \epsilon_j}(-1) e_{\epsilon_i -
\epsilon_j}(-1) {\bf
1} \nonumber \\
&& \quad + \sum_{j=i+1}^{l} e_{\epsilon_1 - \epsilon_j}(-1)
e_{\epsilon_i + \epsilon_j}(-1) {\bf 1} +
\frac{2l-3}{2}e_{\epsilon_1 +
\epsilon_i}(-2) {\bf 1}, \quad \mbox{for } i=2, \ldots ,l, \label{rel.emb.2} \\
&& f_{\epsilon_1}(0)f_{\epsilon_1 - \epsilon_2}(0). v_{B_{l}}=
\frac{1}{2} e_{\epsilon_2}(-1) h_{\epsilon_2}(-1) {\bf 1} -
e_{\epsilon_1 + \epsilon_2}(-1) f_{\epsilon_1}(-1) {\bf 1} +
e_{\epsilon_1}(-1) f_{\epsilon_1 - \epsilon_2}(-1) {\bf 1} \nonumber \\
&& \quad +\sum_{j=3}^{l} e_{\epsilon_2 - \epsilon_j}(-1)
e_{\epsilon_j}(-1) {\bf 1} + \sum_{j=3}^{l} e_{\epsilon_2 +
\epsilon_j}(-1)f_{\epsilon_j}(-1) {\bf 1} - \frac{2l-5}{2}
e_{\epsilon_2}(-2) {\bf 1}, \label{rel.emb.3} \\
&& f_{\epsilon_1}(0)^2. v_{B_{l}}= \frac{1}{2}( e_{\epsilon_1}(-1)
f_{\epsilon_1}(-1) {\bf 1}+ f_{\epsilon_1}(-1) e_{\epsilon_1}(-1)
{\bf 1}) - \frac{1}{2} h_{\epsilon_1}(-1) ^2 {\bf 1} \nonumber \\
&&\quad + \sum_{j=2}^{l}( e_{\epsilon_j}(-1) f_{\epsilon_j}(-1) {\bf
1} + f_{\epsilon_j}(-1) e_{\epsilon_j}(-1) {\bf 1}) \nonumber \\
&& \quad - \sum_{j=2}^{l}( e_{\epsilon_1 + \epsilon_j}(-1)
f_{\epsilon_1 +\epsilon_j}(-1) {\bf 1} +f_{\epsilon_1 +
\epsilon_j}(-1) e_{\epsilon_1 +\epsilon_j}(-1) {\bf 1}) \nonumber \\
&& \quad - \sum_{j=2}^{l}( e_{\epsilon_1 - \epsilon_j}(-1)
f_{\epsilon_1 -\epsilon_j}(-1) {\bf 1} + e_{\epsilon_1 -
\epsilon_j}(-1) f_{\epsilon_1 -\epsilon_j}(-1) {\bf 1}  \label{rel.emb.4} \\
&&f_{\epsilon_1 - \epsilon_i}(0)f_{\epsilon_1 + \epsilon_i}(0).
v_{B_{l}}=  -\frac{1}{4}( e_{\epsilon_i}(-1) f_{\epsilon_i}(-1) {\bf
1}+f_{\epsilon_i}(-1) e_{\epsilon_i}(-1) {\bf 1}) \nonumber \\
&& \quad + \frac{1}{4}( e_{\epsilon_1}(-1) f_{\epsilon_1}(-1) {\bf
1}+ f_{\epsilon_1}(-1) e_{\epsilon_1}(-1) {\bf 1})
 + h_{\epsilon_1 -
\epsilon_i}(-1) h_{\epsilon_1 + \epsilon_i}(-1) {\bf 1} \nonumber \\
&& \quad + \frac{1}{2} \sum_{j=2 \atop j\neq i}^{l} (e_{\epsilon_1 +
\epsilon_j}(-1) f_{\epsilon_1 +\epsilon_j}(-1) {\bf 1}+
f_{\epsilon_1 + \epsilon_j}(-1) e_{\epsilon_1 +\epsilon_j}(-1) {\bf
1}) \nonumber \\
&& \quad + \frac{1}{2}\sum_{j=2 \atop j\neq i}^{l}( e_{\epsilon_1 -
\epsilon_j}(-1) f_{\epsilon_1 -\epsilon_j}(-1) {\bf 1}+
f_{\epsilon_1 - \epsilon_j}(-1) e_{\epsilon_1 -\epsilon_j}(-1) {\bf
1}) \nonumber \\
&& \quad - \frac{1}{2} \sum_{j=2}^{i-1}( e_{\epsilon_j -
\epsilon_i}(-1) f_{\epsilon_j -\epsilon_i}(-1) {\bf 1}
+f_{\epsilon_j - \epsilon_i}(-1) e_{\epsilon_j -\epsilon_i}(-1) {\bf
1}) \nonumber \\
&& \quad - \frac{1}{2} \sum_{j=2 \atop j\neq i}^{l}( e_{\epsilon_i +
\epsilon_j}(-1) f_{\epsilon_i +\epsilon_j}(-1) {\bf 1}
+f_{\epsilon_i + \epsilon_j}(-1) e_{\epsilon_i +\epsilon_j}(-1) {\bf
1}) \nonumber \\
&& \quad  - \frac{1}{2} \sum_{j=i+1}^{l}( e_{\epsilon_i -
\epsilon_j}(-1) f_{\epsilon_i -\epsilon_j}(-1) {\bf 1}
+f_{\epsilon_i - \epsilon_j}(-1) e_{\epsilon_i
-\epsilon_j}(-1) {\bf 1}), \quad \mbox{for } i=2, \ldots ,l, \label{rel.emb.5} \\
 && f_{\epsilon_1 + \epsilon_i}(0). v_{B_{l}}= -\frac{1}{2}
e_{\epsilon_1}(-1) f_{\epsilon_i}(-1) {\bf 1} +\sum_{j=2}^{i-1}
e_{\epsilon_1 + \epsilon_j}(-1) f_{\epsilon_j + \epsilon_i}(-1) {\bf
1} \nonumber \\
&& \quad - \sum_{j=2}^{i-1} e_{\epsilon_1 - \epsilon_j}(-1)
e_{\epsilon_j - \epsilon_i}(-1) {\bf 1} - h_{\epsilon_1 +
\epsilon_i}(-1) e_{\epsilon_1 - \epsilon_i}(-1) {\bf 1}-
\sum_{j=i+1}^{l} e_{\epsilon_1 + \epsilon_j}(-1) f_{\epsilon_i +
\epsilon_j}(-1) {\bf
1} \nonumber \\
&& \quad - \sum_{j=i+1}^{l} e_{\epsilon_1 - \epsilon_j}(-1)
f_{\epsilon_i - \epsilon_j}(-1) {\bf 1} +
\frac{2l-3}{2}e_{\epsilon_1 -
\epsilon_i}(-2) {\bf 1}, \quad \mbox{for } i=3, \ldots ,l, \label{rel.emb.6} \\
&&f_{\epsilon_1 - \epsilon_i}(0)f_{\epsilon_1 - \epsilon_2}(0).
v_{B_{l}}= -\frac{1}{2} e_{\epsilon_2}(-1) e_{\epsilon_i}(-1) {\bf
1}- h_{\epsilon_2 - \epsilon_i}(-1) e_{\epsilon_2 + \epsilon_i}(-1)
{\bf 1} \nonumber \\
&&\quad - f_{\epsilon_1 - \epsilon_2}(-1) e_{\epsilon_1 +
\epsilon_i}(-1) {\bf 1} - e_{\epsilon_1 + \epsilon_2}(-1)
f_{\epsilon_1 - \epsilon_i}(-1) {\bf 1} + \sum_{j=3}^{i-1}
e_{\epsilon_2 + \epsilon_j}(-1) f_{\epsilon_j -\epsilon_i}(-1) {\bf
1} \nonumber \\
&& \quad - \sum_{j=3}^{i-1} e_{\epsilon_2 - \epsilon_j}(-1)
e_{\epsilon_i +\epsilon_j}(-1) {\bf 1} + \sum_{j=i+1}^{l}
e_{\epsilon_2 - \epsilon_j}(-1) e_{\epsilon_i +\epsilon_j}(-1) {\bf
1}  \nonumber \\
&& \quad + \sum_{j=i+1}^{l} e_{\epsilon_2 + \epsilon_j}(-1)
e_{\epsilon_i -\epsilon_j}(-1) {\bf 1} + \frac{2l-3}{2}
e_{\epsilon_2 + \epsilon_i}(-2) {\bf 1}, \quad \mbox{for } i=3, \ldots ,l, \label{rel.emb.7} \\
 &&f_{\epsilon_1 - \epsilon_2}(0)^2. v_{B_{l}}= -\frac{1}{2}
e_{\epsilon_2}(-1)^2 {\bf 1} - 2 f_{\epsilon_1 - \epsilon_2}(-1)
e_{\epsilon_1 + \epsilon_2}(-1) {\bf 1} + 2 \sum_{j=3}^{l}
e_{\epsilon_2 - \epsilon_j}(-1) e_{\epsilon_2 + \epsilon_j}(-1) {\bf
1}, \nonumber \\
&& \mbox{} \label{rel.emb.8} \\
 &&f_{\epsilon_1 - \epsilon_2}(0)f_{\epsilon_1 + \epsilon_i}(0).
v_{B_{l}}= -\frac{1}{2} e_{\epsilon_2}(-1) f_{\epsilon_i}(-1) {\bf
1}- h_{\epsilon_2 + \epsilon_i}(-1) e_{\epsilon_2 - \epsilon_i}(-1)
{\bf 1} \nonumber \\
&&\quad - f_{\epsilon_1 - \epsilon_2}(-1) e_{\epsilon_1 -
\epsilon_i}(-1) {\bf 1} - e_{\epsilon_1 + \epsilon_2}(-1)
f_{\epsilon_1 + \epsilon_i}(-1) {\bf 1} + \sum_{j=3}^{i-1}
e_{\epsilon_2 + \epsilon_j}(-1) f_{\epsilon_i +\epsilon_j}(-1) {\bf
1} \nonumber \\
&& \quad - \sum_{j=3}^{i-1} e_{\epsilon_2 - \epsilon_j}(-1)
e_{\epsilon_j -\epsilon_i}(-1) {\bf 1} - \sum_{j=i+1}^{l}
e_{\epsilon_2 - \epsilon_j}(-1) f_{\epsilon_i -\epsilon_j}(-1) {\bf
1} \nonumber \\
&& \quad - \sum_{j=i+1}^{l} e_{\epsilon_2 + \epsilon_j}(-1)
f_{\epsilon_i +\epsilon_j}(-1) {\bf 1} + \frac{2l-3}{2}
e_{\epsilon_2 - \epsilon_i}(-2) {\bf 1}, \quad \mbox{for } i=3,
\ldots ,l. \label{rel.emb.9}
\end{eqnarray}
Using relations (\ref{rel.emb.1})--(\ref{rel.emb.9}) one can obtain
\begin{eqnarray*}
&& v_{D_{l}}=\Bigg( \frac{2l+1}{12}e_{\epsilon_2}(-1)^2
-\frac{2l-5}{3}f_{\epsilon_1 - \epsilon_2}(-1) e_{\epsilon_1 +
\epsilon_2}(-1)  \\
&&+ \frac{2l+1}{3} \sum_{i=3}^{l} e_{\epsilon_2 - \epsilon_i}(-1)
e_{\epsilon_2 + \epsilon_i}(-1)-\frac{1}{2}e_{\epsilon_1 +
\epsilon_2}(-1)e_{\epsilon_2}(-1) f_{\epsilon_1}(0) \\
&&+ \sum_{i=3}^{l} e_{\epsilon_1 + \epsilon_2}(-1) e_{\epsilon_2 -
\epsilon_i}(-1) f_{\epsilon_1 - \epsilon_i}(0) - \frac{2l+1}{12}
e_{\epsilon_1}(-1)e_{\epsilon_2}(-1)f_{\epsilon_1 - \epsilon_2}(0)
\\
&& - \frac{(2l+1)(2l-5)}{12}e_{\epsilon_1 +
\epsilon_2}(-2)f_{\epsilon_1 - \epsilon_2}(0)+ \frac{2l-5}{6}
h_{\epsilon_1 - \epsilon_2}(-1)e_{\epsilon_1 +
\epsilon_2}(-1)f_{\epsilon_1 - \epsilon_2}(0)
\\
&& - \frac{2l+1}{6} \sum_{i=3}^{l} e_{\epsilon_1 - \epsilon_i}(-1)
e_{\epsilon_2 + \epsilon_i}(-1) f_{\epsilon_1 -
\epsilon_2}(0) \\
&& - \frac{2l+1}{6} \sum_{i=3}^{l} e_{\epsilon_1 + \epsilon_i}(-1)
e_{\epsilon_2 - \epsilon_i}(-1) f_{\epsilon_1 - \epsilon_2}(0) +
\frac{1}{2} e_{\epsilon_1}(-1) e_{\epsilon_1 +
\epsilon_2}(-1)f_{\epsilon_1}(0) f_{\epsilon_1 - \epsilon_2}(0)\\
&& +\frac{1}{2l-1}e_{\epsilon_1 + \epsilon_2}(-1)^2
f_{\epsilon_1}(0)^2 -\frac{2l-5}{2l-1} e_{\epsilon_1 +
\epsilon_2}(-1)^2 f_{\epsilon_1 - \epsilon_2}(0) f_{\epsilon_1 +
\epsilon_2}(0)\\
&& + \frac{4}{2l-1} \sum_{i=3}^{l} e_{\epsilon_1 + \epsilon_2}(-1)^2
f_{\epsilon_1 - \epsilon_i}(0) f_{\epsilon_1 + \epsilon_i}(0) \\
&& + \sum_{i=3}^{l} e_{\epsilon_1 + \epsilon_2}(-1) e_{\epsilon_2 +
\epsilon_i}(-1) f_{\epsilon_1 + \epsilon_i}(0)- \sum_{i=3}^{l}
e_{\epsilon_1 + \epsilon_2}(-1) e_{\epsilon_1 - \epsilon_i}(-1)
f_{\epsilon_1 - \epsilon_i}(0) f_{\epsilon_1 - \epsilon_2}(0) \\
&&+ \frac{2l+1}{24} e_{\epsilon_1}(-1)^2 f_{\epsilon_1 -
\epsilon_2}(0)^2 + \frac{2l-5}{6} e_{\epsilon_1 - \epsilon_2}(-1)
e_{\epsilon_1 + \epsilon_2}(-1)f_{\epsilon_1 - \epsilon_2}(0)^2 \\
&& + \frac{2l+1}{6} \sum_{i=3}^{l} e_{\epsilon_1 - \epsilon_i}(-1)
e_{\epsilon_1 + \epsilon_i}(-1) f_{\epsilon_1 -
\epsilon_2}(0)^2 \\
&& - \sum_{i=3}^{l} e_{\epsilon_1 + \epsilon_2}(-1) e_{\epsilon_1 +
\epsilon_i}(-1) f_{\epsilon_1 - \epsilon_2}(0) f_{\epsilon_1 +
\epsilon_i}(0) \Bigg). v_{B_{l}},
\end{eqnarray*}
which implies the claim of lemma. $\;\;\;\;\Box$

It follows from Lemma \ref{l.emb.tech}, Theorem \ref{t.3.6} and
Proposition \ref{p.3.1} that the embedding of
$N_{D_{l}}(-l+\frac{3}{2},0)$ into $N_{B_{l}}(-l+\frac{3}{2},0)$
induces the embedding of $L_{D_{l}}(-l+\frac{3}{2},0)$ into
$L_{B_{l}}(-l+\frac{3}{2},0)$. We get:

\begin{thm}
$L_{D_{l}}(-l+\frac{3}{2},0)$ is a vertex subalgebra of
$L_{B_{l}}(-l+\frac{3}{2},0)$.
\end{thm}

Next, we show that $L_{D_{l}}(-l+\frac{3}{2},0)$ is a vertex
operator subalgebra of $L_{B_{l}}(-l+\frac{3}{2},0)$, i.e. that
these vertex operator algebras have the same conformal vector. We
use the following lemma:
\begin{lem} \label{lem.dokaz.Vir}
 Relation
\begin{eqnarray*}
&& (2l-1) \sum_{\alpha \in \Delta _{B_{l}}^{+} \atop
(\alpha,\alpha)=1 }( e_{\alpha}(-1)f_{\alpha}(-1)
{\bf 1}+f_{\alpha}(-1)e_{\alpha}(-1) {\bf 1})  \\
&& =4\sum_{\alpha \in \Delta _{D_{l}}^{+} }( e_{\alpha}(-1)
f_{\alpha}(-1){\bf 1} +f_{\alpha}(-1)e_{\alpha}(-1) {\bf 1}) +
\sum_{\alpha \in \Delta _{B_{l}}^{+} \atop (\alpha,\alpha)=1 }
h_{\alpha}(-1)^{2}{\bf 1} \qquad
\end{eqnarray*}
holds in $L_{B_{l}}(-l +\frac{3}{2},0)$.
\end{lem}
{\bf Proof:} Proposition \ref{p.3.1} implies that relation
\begin{eqnarray*}
&& \Big( 2l f_{\epsilon_1}(0)^2 + 4 \sum _{i=2}^{l}f_{\epsilon_1
-\epsilon_i}(0) f_{\epsilon_1 + \epsilon_i}(0) \Big). v_{B_{l}} =0
\nonumber
\end{eqnarray*}
holds in $L_{B_{l}}(-l+\frac{3}{2},0)$. Using relations
(\ref{rel.emb.4}) and (\ref{rel.emb.5}) we get the claim of lemma.
$\;\;\;\;\Box$

\begin{thm} \label{prop.jedn.Vir}
Denote by $\omega _{B_{l}}$ the conformal vector for vertex operator
algebra $L_{B_{l}}(-l+\frac{3}{2},0)$, and by $\omega _{D_{l}}$ the
conformal vector for $L_{D_{l}}(-l+\frac{3}{2},0)$. Then
\begin{eqnarray*}
\omega _{B_{l}}=\omega _{D_{l}}.
\end{eqnarray*}
\end{thm}
{\bf Proof:} Set $\{ \frac{1}{2} h_{\alpha} \ \vert \ \alpha \in
\Delta _{B_{l}}^{+}, (\alpha,\alpha)=1 \}$ is an orthonormal basis
of ${\frak h}$ with respect to the form $(\cdot, \cdot)$. Lemma
\ref{lem.dokaz.Vir} implies that
\begin{eqnarray} \label{rel.dokaz.Vir4}
&& \sum_{\alpha \in \Delta _{B_{l}}^{+} \atop (\alpha,\alpha)=1 }(
e_{\alpha}(-1)f_{\alpha}(-1)
{\bf 1}+f_{\alpha}(-1)e_{\alpha}(-1) {\bf 1}) \nonumber \\
&& = \frac{4}{2l-1} \sum_{\alpha \in \Delta _{D_{l}}^{+} }(
e_{\alpha}(-1) f_{\alpha}(-1){\bf 1} +f_{\alpha}(-1)e_{\alpha}(-1)
{\bf 1}) + \frac{1}{2l-1} \sum_{\alpha \in \Delta _{B_{l}}^{+} \atop
(\alpha,\alpha)=1 } h_{\alpha}(-1)^{2}{\bf 1} \nonumber \\
&& \mbox{}
\end{eqnarray}
holds in $L_{B_{l}}(-l +\frac{3}{2},0)$. It follows from relation
(\ref{rel.Virasoro}) that
\begin{eqnarray*}
&& \omega _{B_{l}}=\frac{1}{2l+1} \Bigg(\frac{1}{4} \sum_{\alpha \in
\Delta _{B_{l}}^{+} \atop (\alpha,\alpha)=1 } h_{\alpha}(-1)^{2}{\bf
1} + \sum_{\alpha \in \Delta _{B_{l}}^{+}} \frac{(\alpha,
\alpha)}{2}(e_{\alpha}(-1)f_{\alpha}(-1) {\bf 1}+
f_{\alpha}(-1)e_{\alpha}(-1) {\bf 1}) \Bigg) \\
&& = \frac{1}{2l+1} \Bigg( \frac{1}{4} \sum_{\alpha \in \Delta
_{B_{l}}^{+} \atop (\alpha,\alpha)=1 } h_{\alpha}(-1)^{2}{\bf 1}+
\sum_{\alpha \in \Delta _{D_{l}}^{+} }
(e_{\alpha}(-1)f_{\alpha}(-1) {\bf 1}+ f_{\alpha}(-1)e_{\alpha}(-1) {\bf 1}) \\
&&  \qquad + \frac{1}{2} \sum_{\alpha \in \Delta _{B_{l}}^{+} \atop
(\alpha,\alpha)=1 }( e_{\alpha}(-1)f_{\alpha}(-1)+
f_{\alpha}(-1)e_{\alpha}(-1) {\bf 1}) \Bigg).
\end{eqnarray*}
Using relation (\ref{rel.dokaz.Vir4}), we obtain
\begin{eqnarray*}
&& \omega _{B_{l}}=\frac{1}{2l-1} \Bigg(\frac{1}{4} \sum_{\alpha \in
\Delta _{B_{l}}^{+} \atop (\alpha,\alpha)=1 } h_{\alpha}(-1)^{2}{\bf
1} + \sum_{\alpha \in \Delta _{D_{l}}^{+}}
(e_{\alpha}(-1)f_{\alpha}(-1) {\bf 1}+ f_{\alpha}(-1)e_{\alpha}(-1)
{\bf 1}) \Bigg)= \omega _{D_{l}}. \;\;\;\;\Box
\end{eqnarray*}

Thus, $L_{D_{l}}(-l+\frac{3}{2},0)$ is a vertex operator subalgebra
of $L_{B_{l}}(-l+\frac{3}{2},0)$.

\section{Case $l=4$}

For $l=4$, simple Lie algebra $\frak g _{F_{4}}$ contains three
copies of $\frak g _{B_{4}}$ as Lie subalgebras. We denote these
copies $\frak g _{B_{4}}$, $\frak g _{B_{4}}'$ and $\frak g
_{B_{4}}''$, and we have the corresponding isomorphisms
\begin{eqnarray*}
&& \pi ' : \frak g _{B_{4}} \to \frak g _{B_{4}}',
\\
&& \pi '' : \frak g _{B_{4}} \to \frak g _{B_{4}}''.
\end{eqnarray*}
Lie algebra $\frak g _{D_{4}}$ is a Lie subalgebra of all three
copies of $\frak g _{B_{4}}$, and it is invariant under isomorphisms
$\pi '$ and $\pi ''$. When restricted to $\frak g _{D_{4}}$, $\pi '$
and $\pi ''$ are automorphisms induced from the automorphisms of the
root system $\Delta _{D_{4}}$, determined by:
\begin{eqnarray}
&& \pi '(\alpha_1)=\alpha_3, \  \pi '(\alpha_2)=\alpha_2, \ \pi
'(\alpha_3)=\alpha_4, \ \pi '(\alpha_4)=\alpha_1 \nonumber \\
&& \pi ''(\alpha_1)=\alpha_4, \  \pi ''(\alpha_2)=\alpha_2, \ \pi
''(\alpha_3)=\alpha_3, \ \pi ''(\alpha_4)=\alpha_1. \label{rel.pi}
\end{eqnarray}

Vertex operator algebra $L_{F_{4}}(-\frac{5}{2},0)$ contains three
copies of $L_{B_{4}}(-\frac{5}{2},0)$ as vertex operator subalgebras
(cf. \cite{P2}). We denote these copies $L_{B_{4}}(-\frac{5}{2},0)$,
$L_{B_{4}}'(-\frac{5}{2},0)$ and $L_{B_{4}}''(-\frac{5}{2},0)$, and
the corresponding isomorphisms of vertex operator algebras are
induced from $\pi '$ and $\pi ''$. Using relations (\ref{rel.pi})
and (\ref{rel.sing.v.D}), one can easily verify that $\pi '
(v_{D_{4}})=v_{D_{4}}$ and $\pi '' (v_{D_{4}})=v_{D_{4}}$, which
implies that $L_{D_{4}}(-\frac{5}{2},0)$ is a vertex subalgebra of
all three copies $L_{B_{4}}(-\frac{5}{2},0)$,
$L_{B_{4}}'(-\frac{5}{2},0)$ and $L_{B_{4}}''(-\frac{5}{2},0)$. We
obtain:

\begin{coro}
$L_{D_{4}}(-\frac{5}{2},0)$ is a vertex subalgebra of three copies
of $L_{B_{4}}(-\frac{5}{2},0)$ contained in
$L_{F_{4}}(-\frac{5}{2},0)$. Vertex operator algebras
$L_{D_{4}}(-\frac{5}{2},0)$, $L_{B_{4}}(-\frac{5}{2},0)$,
$L_{B_{4}}'(-\frac{5}{2},0)$, $L_{B_{4}}''(-\frac{5}{2},0)$ and
$L_{F_{4}}(-\frac{5}{2},0)$ have the same conformal vector.
\end{coro}

\section{Appendix}

In this Appendix we show relation $e_{\epsilon_1 -
\epsilon_{2}}(0).v_{D_{l}}=0$ from the proof of Lemma \ref{l.3.5}.
Relations $e_{\epsilon_2 - \epsilon_{3}}(0).v_{D_{l}}=0$,
$e_{\epsilon_k - \epsilon_{k+1}}(0).v_{D_{l}}=0$, for $k=3,\ldots
,l-1$, $e_{\epsilon_{l-1} + \epsilon_{l}}(0).v_{D_{l}}=0$ and
$f_{\theta}(1).v_{D_{l}}=0$  can be verified similarly. We have:
\begin{eqnarray}
&& e_{\epsilon_1 - \epsilon_{2}}(0). e_{\epsilon_1 - \epsilon_i}(-1)
e_{\epsilon_1 + \epsilon_i}(-1) e_{\epsilon_2 - \epsilon_j}(-1)
e_{\epsilon_2 + \epsilon_j}(-1) {\bf 1} \nonumber \\
&& \ = e_{\epsilon_1 - \epsilon_i}(-1) e_{\epsilon_1 +
\epsilon_i}(-1) e_{\epsilon_1 - \epsilon_j}(-1) e_{\epsilon_2 +
\epsilon_j}(-1) {\bf 1} \nonumber \\
&& \ +e_{\epsilon_1 - \epsilon_i}(-1) e_{\epsilon_1 +
\epsilon_i}(-1) e_{\epsilon_1 + \epsilon_j}(-1) e_{\epsilon_2 -
\epsilon_j}(-1) {\bf 1} \nonumber \\
&& \ + e_{\epsilon_1 - \epsilon_i}(-1) e_{\epsilon_1 +
\epsilon_i}(-1)  e_{\epsilon_1 + \epsilon_2}(-2) {\bf 1}, \quad
\mbox{for } i,j=3, \ldots ,l, \label{rel.app.1}\\
&& e_{\epsilon_1 - \epsilon_{2}}(0). e_{\epsilon_1 - \epsilon_i}(-1)
e_{\epsilon_2 + \epsilon_i}(-1) e_{\epsilon_1 + \epsilon_j}(-1)
e_{\epsilon_2 - \epsilon_j}(-1) {\bf 1} \nonumber \\
&& \ = e_{\epsilon_1 - \epsilon_i}(-1) e_{\epsilon_1 +
\epsilon_i}(-1) e_{\epsilon_1 + \epsilon_j}(-1) e_{\epsilon_2 -
\epsilon_j}(-1) {\bf 1} \nonumber \\
&& \ + e_{\epsilon_1 - \epsilon_i}(-1) e_{\epsilon_2 +
\epsilon_i}(-1) e_{\epsilon_1 + \epsilon_j}(-1) e_{\epsilon_1 -
\epsilon_j}(-1) {\bf 1}, \quad \mbox{for } i,j=3, \ldots ,l,  i \neq
j, \nonumber \\
&& \mbox{} \\
&& e_{\epsilon_1 - \epsilon_{2}}(0). e_{\epsilon_1 - \epsilon_i}(-1)
e_{\epsilon_2 + \epsilon_i}(-1) e_{\epsilon_1 - \epsilon_j}(-1)
e_{\epsilon_2 + \epsilon_j}(-1) {\bf 1} \nonumber \\
&& \ = e_{\epsilon_1 - \epsilon_i}(-1) e_{\epsilon_1 +
\epsilon_i}(-1) e_{\epsilon_1 - \epsilon_j}(-1) e_{\epsilon_2 +
\epsilon_j}(-1) {\bf 1} \nonumber \\
&& \ + e_{\epsilon_1 - \epsilon_i}(-1) e_{\epsilon_2 +
\epsilon_i}(-1) e_{\epsilon_1 - \epsilon_j}(-1) e_{\epsilon_1 +
\epsilon_j}(-1) {\bf 1}, \quad \mbox{for } i,j=3, \ldots ,l,  i \neq
j, \nonumber \\
&& \mbox{} \\
&& e_{\epsilon_1 - \epsilon_{2}}(0). e_{\epsilon_1 - \epsilon_i}(-1)
e_{\epsilon_2 + \epsilon_i}(-1) e_{\epsilon_1 - \epsilon_i}(-1)
e_{\epsilon_2 + \epsilon_i}(-1) {\bf 1} \nonumber \\
&& \ = 2 e_{\epsilon_1 - \epsilon_i}(-1)^2 e_{\epsilon_1 +
\epsilon_i}(-1) e_{\epsilon_2 + \epsilon_i}(-1)
 {\bf 1} + e_{\epsilon_1 + \epsilon_2}(-2) e_{\epsilon_1 -
\epsilon_i}(-1) e_{\epsilon_1 + \epsilon_i}(-1)  {\bf 1}  \nonumber \\
&& \quad \mbox{for } i=3, \ldots ,l, \\
&& e_{\epsilon_1 - \epsilon_{2}}(0). e_{\epsilon_1 + \epsilon_i}(-1)
e_{\epsilon_2 - \epsilon_i}(-1) e_{\epsilon_1 + \epsilon_j}(-1)
e_{\epsilon_2 - \epsilon_j}(-1) {\bf 1} \nonumber \\
&& \ = e_{\epsilon_1 + \epsilon_i}(-1) e_{\epsilon_1 -
\epsilon_i}(-1) e_{\epsilon_1 + \epsilon_j}(-1) e_{\epsilon_2 -
\epsilon_j}(-1) {\bf 1} \nonumber \\
&& \ + e_{\epsilon_1 + \epsilon_i}(-1) e_{\epsilon_2 -
\epsilon_i}(-1) e_{\epsilon_1 + \epsilon_j}(-1) e_{\epsilon_1 -
\epsilon_j}(-1) {\bf 1}, \quad \mbox{for } i,j=3, \ldots ,l,  i \neq
j, \nonumber
\\
&& \mbox{} \\
&& e_{\epsilon_1 - \epsilon_{2}}(0). e_{\epsilon_1 + \epsilon_i}(-1)
e_{\epsilon_2 - \epsilon_i}(-1) e_{\epsilon_1 + \epsilon_i}(-1)
e_{\epsilon_2 - \epsilon_i}(-1) {\bf 1} \nonumber \\
&& \ = 2 e_{\epsilon_1 + \epsilon_i}(-1)^2 e_{\epsilon_1 -
\epsilon_i}(-1) e_{\epsilon_2 - \epsilon_i}(-1)
 {\bf 1} + e_{\epsilon_1 + \epsilon_2}(-2) e_{\epsilon_1 +
\epsilon_i}(-1) e_{\epsilon_1 - \epsilon_i}(-1)  {\bf 1}  \nonumber \\
&& \quad \mbox{for } i=3, \ldots ,l, \\
&& e_{\epsilon_1 - \epsilon_{2}}(0). e_{\epsilon_1 + \epsilon_2}(-1)
e_{\epsilon_1 + \epsilon_j}(-1) e_{\epsilon_2 - \epsilon_i}(-1)
f_{\epsilon_j - \epsilon_i}(-1) {\bf 1} \nonumber \\
&& \ = e_{\epsilon_1 + \epsilon_2}(-1) e_{\epsilon_1 +
\epsilon_j}(-1) e_{\epsilon_1 - \epsilon_i}(-1) f_{\epsilon_j -
\epsilon_i}(-1) {\bf
1}, \quad \mbox{for } i,j=3, \ldots ,l,  j < i, \nonumber \\
&& \mbox{} \\
&& e_{\epsilon_1 - \epsilon_{2}}(0).e_{\epsilon_1 + \epsilon_2}(-1)
e_{\epsilon_1 + \epsilon_j}(-1) e_{\epsilon_2 - \epsilon_i}(-1)
e_{\epsilon_i - \epsilon_j}(-1) {\bf 1} \nonumber \\
&& \ = e_{\epsilon_1 + \epsilon_2}(-1) e_{\epsilon_1 +
\epsilon_j}(-1) e_{\epsilon_1 - \epsilon_i}(-1) e_{\epsilon_i -
\epsilon_j}(-1) {\bf 1} \quad \mbox{for } i,j=3, \ldots ,l,  i < j, \nonumber \\
&& \mbox{} \\
&& e_{\epsilon_1 - \epsilon_{2}}(0). e_{\epsilon_1 + \epsilon_2}(-1)
e_{\epsilon_1 - \epsilon_j}(-1) e_{\epsilon_2 -
\epsilon_i}(-1) e_{\epsilon_i + \epsilon_j}(-1) {\bf 1} \nonumber \\
&& \ = e_{\epsilon_1 + \epsilon_2}(-1) e_{\epsilon_1 -
\epsilon_j}(-1) e_{\epsilon_1 - \epsilon_i}(-1) e_{\epsilon_i +
\epsilon_j}(-1) {\bf 1}, \quad \mbox{for } i,j=3, \ldots ,l,  i \neq
j, \nonumber
\\
&& \mbox{} \\
&& e_{\epsilon_1 - \epsilon_{2}}(0).e_{\epsilon_1 + \epsilon_2}(-1)
e_{\epsilon_2 + \epsilon_i}(-1) e_{\epsilon_1 +
\epsilon_j}(-1) f_{\epsilon_i + \epsilon_j}(-1) {\bf 1} \nonumber \\
&& \ = e_{\epsilon_1 + \epsilon_2}(-1) e_{\epsilon_1 +
\epsilon_i}(-1) e_{\epsilon_1 + \epsilon_j}(-1) f_{\epsilon_i +
\epsilon_j}(-1) {\bf 1}, \quad \mbox{for } i,j=3, \ldots ,l,  i \neq
j, \nonumber
\\
&& \mbox{} \\
&& e_{\epsilon_1 - \epsilon_{2}}(0). e_{\epsilon_1 + \epsilon_2}(-1)
e_{\epsilon_2 + \epsilon_i}(-1) e_{\epsilon_1 - \epsilon_j}(-1)
e_{\epsilon_j - \epsilon_i}(-1) {\bf 1} \nonumber \\
&& \ = e_{\epsilon_1 + \epsilon_2}(-1) e_{\epsilon_1 +
\epsilon_i}(-1) e_{\epsilon_1 - \epsilon_j}(-1) e_{\epsilon_j -
\epsilon_i}(-1) {\bf 1}, \quad \mbox{for } i,j=3, \ldots ,l,  j<i,
\nonumber
\\
&& \mbox{} \\
&& e_{\epsilon_1 - \epsilon_{2}}(0). e_{\epsilon_1 + \epsilon_2}(-1)
e_{\epsilon_2 + \epsilon_i}(-1) e_{\epsilon_1 -
\epsilon_j}(-1) f_{\epsilon_i - \epsilon_j}(-1) {\bf 1} \nonumber \\
&& \ = e_{\epsilon_1 + \epsilon_2}(-1) e_{\epsilon_1 +
\epsilon_i}(-1) e_{\epsilon_1 - \epsilon_j}(-1) f_{\epsilon_i -
\epsilon_j}(-1) {\bf 1}, \quad \mbox{for } i,j=3, \ldots ,l,  i<j,
\nonumber
\\
&& \mbox{} \\
&& e_{\epsilon_1 - \epsilon_{2}}(0).e_{\epsilon_1 + \epsilon_2}(-1)
f_{\epsilon_1 - \epsilon_2}(-1) e_{\epsilon_1 - \epsilon_i}(-1)
e_{\epsilon_1 +
\epsilon_i}(-1) {\bf 1} \nonumber \\
&& \ = e_{\epsilon_1 + \epsilon_2}(-1)  e_{\epsilon_1 +
\epsilon_i}(-1) e_{\epsilon_1 - \epsilon_i}(-2) {\bf 1} +
e_{\epsilon_1 + \epsilon_2}(-1)  e_{\epsilon_1 - \epsilon_i}(-1)
e_{\epsilon_1 + \epsilon_i}(-2) {\bf 1} \nonumber \\
&& \ + e_{\epsilon_1 + \epsilon_2}(-1)  e_{\epsilon_1 -
\epsilon_i}(-1) e_{\epsilon_1 + \epsilon_i}(-1) h_{\epsilon_1 -
\epsilon_2}(-1) {\bf 1}, \quad \mbox{for } i=3, \ldots ,l, \\
&& e_{\epsilon_1 - \epsilon_{2}}(0). e_{\epsilon_1 + \epsilon_2}(-1)
e_{\epsilon_1 - \epsilon_2}(-1) e_{\epsilon_2 - \epsilon_i}(-1)
e_{\epsilon_2 +
\epsilon_i}(-1) {\bf 1} \nonumber \\
&& \ = e_{\epsilon_1 + \epsilon_2}(-1) e_{\epsilon_1 -
\epsilon_2}(-1) e_{\epsilon_1 - \epsilon_i}(-1) e_{\epsilon_2 +
\epsilon_i}(-1) {\bf 1} \nonumber \\
&& \ + e_{\epsilon_1 + \epsilon_2}(-1) e_{\epsilon_1 -
\epsilon_2}(-1) e_{\epsilon_1 + \epsilon_2}(-2) {\bf 1} \nonumber \\
&& \ + e_{\epsilon_1 + \epsilon_2}(-1) e_{\epsilon_1 -
\epsilon_2}(-1) e_{\epsilon_1 + \epsilon_i}(-1) e_{\epsilon_2 -
\epsilon_i}(-1) {\bf 1}, \quad \mbox{for } i=3, \ldots ,l, \\
&& e_{\epsilon_1 - \epsilon_{2}}(0). e_{\epsilon_1 + \epsilon_2}(-1)
e_{\epsilon_1 + \epsilon_i}(-1) e_{\epsilon_2 - \epsilon_i}(-1)
h_{\epsilon_i}(-1) {\bf 1} \nonumber \\
&& \ = e_{\epsilon_1 + \epsilon_2}(-1) e_{\epsilon_1 +
\epsilon_i}(-1) e_{\epsilon_1 - \epsilon_i}(-1) h_{\epsilon_i}(-1)
{\bf 1}, \quad \mbox{for } i=3, \ldots ,l, \\
&& e_{\epsilon_1 - \epsilon_{2}}(0). e_{\epsilon_1 + \epsilon_2}(-1)
e_{\epsilon_1 + \epsilon_i}(-1) e_{\epsilon_2 - \epsilon_i}(-1)
h_{\epsilon_1 - \epsilon_2}(-1) {\bf 1} \nonumber \\
&& \ = e_{\epsilon_1 + \epsilon_2}(-1) e_{\epsilon_1 +
\epsilon_i}(-1) e_{\epsilon_1 - \epsilon_i}(-1) h_{\epsilon_1 -
\epsilon_2}(-1) {\bf 1} \nonumber \\
&& \ -2 e_{\epsilon_1 + \epsilon_2}(-1) e_{\epsilon_1 -
\epsilon_2}(-1) e_{\epsilon_1 + \epsilon_i}(-1) e_{\epsilon_2 -
\epsilon_i}(-1) {\bf 1} \nonumber \\
&& \ +2 e_{\epsilon_1 + \epsilon_2}(-1) e_{\epsilon_1 +
\epsilon_i}(-1) e_{\epsilon_1 - \epsilon_i}(-2) {\bf 1},
\quad \mbox{for } i=3, \ldots ,l, \\
&& e_{\epsilon_1 - \epsilon_{2}}(0).e_{\epsilon_1 + \epsilon_2}(-1)
e_{\epsilon_1 - \epsilon_i}(-1) e_{\epsilon_2 + \epsilon_i}(-1)
h_{\epsilon_i}(-1) {\bf 1} \nonumber \\
&& \ = e_{\epsilon_1 + \epsilon_2}(-1) e_{\epsilon_1 -
\epsilon_i}(-1) e_{\epsilon_1 + \epsilon_i}(-1) h_{\epsilon_i}(-1)
{\bf 1}, \quad \mbox{for } i=3, \ldots ,l, \\
&& e_{\epsilon_1 - \epsilon_{2}}(0). e_{\epsilon_1 + \epsilon_2}(-1)
e_{\epsilon_1 - \epsilon_i}(-1) e_{\epsilon_2 + \epsilon_i}(-1)
h_{\epsilon_1 - \epsilon_2}(-1) {\bf 1} \nonumber \\
&& \ = e_{\epsilon_1 + \epsilon_2}(-1) e_{\epsilon_1 -
\epsilon_i}(-1) e_{\epsilon_1 + \epsilon_i}(-1) h_{\epsilon_1 -
\epsilon_2}(-1) {\bf 1} \nonumber \\
&& \ -2 e_{\epsilon_1 + \epsilon_2}(-1) e_{\epsilon_1 -
\epsilon_2}(-1) e_{\epsilon_1 - \epsilon_i}(-1) e_{\epsilon_2 +
\epsilon_i}(-1) {\bf 1} \nonumber \\
&& \ +2 e_{\epsilon_1 + \epsilon_2}(-1) e_{\epsilon_1 -
\epsilon_i}(-1) e_{\epsilon_1 + \epsilon_i}(-2) {\bf 1},
\quad \mbox{for } i=3, \ldots ,l, \\
&& e_{\epsilon_1 - \epsilon_{2}}(0). e_{\epsilon_1 + \epsilon_2}(-1)
e_{\epsilon_1 + \epsilon_i}(-2) e_{\epsilon_2 - \epsilon_i}(-1) {\bf
1} = e_{\epsilon_1 + \epsilon_2}(-1) e_{\epsilon_1 + \epsilon_i}(-2)
e_{\epsilon_1 - \epsilon_i}(-1) {\bf 1} \nonumber \\
&& \quad \mbox{for } i=3, \ldots ,l, \\
&& e_{\epsilon_1 - \epsilon_{2}}(0). e_{\epsilon_1 + \epsilon_2}(-1)
e_{\epsilon_1 - \epsilon_i}(-2) e_{\epsilon_2 + \epsilon_i}(-1) {\bf
1} = e_{\epsilon_1 + \epsilon_2}(-1) e_{\epsilon_1 - \epsilon_i}(-2)
e_{\epsilon_1 + \epsilon_i}(-1) {\bf 1} \nonumber \\
&& \quad \mbox{for } i=3, \ldots ,l, \\
&& e_{\epsilon_1 - \epsilon_{2}}(0). e_{\epsilon_1 + \epsilon_2}(-2)
e_{\epsilon_1 + \epsilon_i}(-1) e_{\epsilon_2 - \epsilon_i}(-1) {\bf
1}= e_{\epsilon_1 + \epsilon_2}(-2) e_{\epsilon_1 + \epsilon_i}(-1)
e_{\epsilon_1 - \epsilon_i}(-1) {\bf 1} \nonumber \\
&& \quad \mbox{for } i=3, \ldots ,l, \\
&& e_{\epsilon_1 - \epsilon_{2}}(0). e_{\epsilon_1 + \epsilon_2}(-2)
e_{\epsilon_1 - \epsilon_i}(-1) e_{\epsilon_2 + \epsilon_i}(-1) {\bf
1}= e_{\epsilon_1 + \epsilon_2}(-2) e_{\epsilon_1 - \epsilon_i}(-1)
e_{\epsilon_1 + \epsilon_i}(-1) {\bf 1} \nonumber \\
&& \quad \mbox{for } i=3, \ldots ,l, \\
&& e_{\epsilon_1 - \epsilon_{2}}(0). e_{\epsilon_1 + \epsilon_2}(-1)
e_{\epsilon_1 + \epsilon_i}(-1) e_{\epsilon_2 - \epsilon_i}(-2) {\bf
1}= e_{\epsilon_1 + \epsilon_2}(-1) e_{\epsilon_1 + \epsilon_i}(-1)
e_{\epsilon_1 - \epsilon_i}(-2) {\bf 1} \nonumber \\
&& \quad \mbox{for } i=3, \ldots ,l, \\
&& e_{\epsilon_1 - \epsilon_{2}}(0). e_{\epsilon_1 + \epsilon_2}(-1)
e_{\epsilon_1 - \epsilon_i}(-1) e_{\epsilon_2 + \epsilon_i}(-2) {\bf
1}= e_{\epsilon_1 + \epsilon_2}(-1) e_{\epsilon_1 - \epsilon_i}(-1)
e_{\epsilon_1 + \epsilon_i}(-2) {\bf 1} \nonumber \\
&& \quad \mbox{for } i=3, \ldots ,l, \\
&& e_{\epsilon_1 - \epsilon_{2}}(0). e_{\epsilon_1 + \epsilon_2}(-2)
e_{\epsilon_1 + \epsilon_2}(-1) h_{\epsilon_1 - \epsilon_2}(-1) {\bf
1}=-2 e_{\epsilon_1 + \epsilon_2}(-2) e_{\epsilon_1 +
\epsilon_2}(-1) e_{\epsilon_1 - \epsilon_2}(-1) {\bf 1}, \nonumber \\
&& \mbox{} \\
&& e_{\epsilon_1 - \epsilon_{2}}(0). e_{\epsilon_1 + \epsilon_2}(-2)
e_{\epsilon_1 + \epsilon_2}(-1) h_{\epsilon_1}(-1) {\bf 1}=-2
e_{\epsilon_1 + \epsilon_2}(-2) e_{\epsilon_1 +
\epsilon_2}(-1) e_{\epsilon_1 - \epsilon_2}(-1) {\bf 1}, \nonumber \\
&& \mbox{} \\
&& e_{\epsilon_1 - \epsilon_{2}}(0). e_{\epsilon_1 +
\epsilon_2}(-2)^2 {\bf 1}=0, \\
&& e_{\epsilon_1 - \epsilon_{2}}(0). e_{\epsilon_1 + \epsilon_2}(-1)
e_{\epsilon_1 + \epsilon_2}(-3) {\bf 1}=0, \\
&& e_{\epsilon_1 - \epsilon_{2}}(0). e_{\epsilon_1 +
\epsilon_2}(-1)^2 (e_{\epsilon_1 - \epsilon_2}(-1)f_{\epsilon_1 -
\epsilon_2}(-1)+ f_{\epsilon_1 - \epsilon_2}(-1)e_{\epsilon_1 -
\epsilon_2}(-1)){\bf 1} \nonumber \\
&& \ = 2e_{\epsilon_1 + \epsilon_2}(-1)^2 e_{\epsilon_1 -
\epsilon_2}(-1)h_{\epsilon_1 - \epsilon_2}(-1) {\bf 1} +2
e_{\epsilon_1 + \epsilon_2}(-1)^2 e_{\epsilon_1 - \epsilon_2}(-2)
{\bf 1}, \\
&& e_{\epsilon_1 - \epsilon_{2}}(0). e_{\epsilon_1 +
\epsilon_2}(-1)^2 (e_{\epsilon_2 - \epsilon_i}(-1)f_{\epsilon_2 -
\epsilon_i}(-1)+ f_{\epsilon_2 - \epsilon_i}(-1)e_{\epsilon_2 -
\epsilon_i}(-1)){\bf 1} \nonumber \\
&& \ = 2e_{\epsilon_1 + \epsilon_2}(-1)^2 e_{\epsilon_1 -
\epsilon_i}(-1)f_{\epsilon_2 - \epsilon_i}(-1) {\bf 1} -
e_{\epsilon_1 + \epsilon_2}(-1)^2 e_{\epsilon_1 - \epsilon_2}(-2)
{\bf 1}, \nonumber \\
&& \quad \mbox{for } i=3, \ldots ,l, \\
&& e_{\epsilon_1 - \epsilon_{2}}(0). e_{\epsilon_1 +
\epsilon_2}(-1)^2 (e_{\epsilon_2 + \epsilon_i}(-1)f_{\epsilon_2 +
\epsilon_i}(-1)+ f_{\epsilon_2 + \epsilon_i}(-1)e_{\epsilon_2 +
\epsilon_i}(-1)){\bf 1} \nonumber \\
&& \ = 2e_{\epsilon_1 + \epsilon_2}(-1)^2 e_{\epsilon_1 +
\epsilon_i}(-1)f_{\epsilon_2 + \epsilon_i}(-1) {\bf 1} -
e_{\epsilon_1 + \epsilon_2}(-1)^2 e_{\epsilon_1 - \epsilon_2}(-2)
{\bf 1}, \nonumber \\
&& \quad \mbox{for } i=3, \ldots ,l, \\
&& e_{\epsilon_1 - \epsilon_{2}}(0). e_{\epsilon_1 +
\epsilon_2}(-1)^2 (e_{\epsilon_1 - \epsilon_i}(-1)f_{\epsilon_1 -
\epsilon_i}(-1)+ f_{\epsilon_1 - \epsilon_i}(-1)e_{\epsilon_1 -
\epsilon_i}(-1)){\bf 1} \nonumber \\
&& \ = -2e_{\epsilon_1 + \epsilon_2}(-1)^2 e_{\epsilon_1 -
\epsilon_i}(-1)f_{\epsilon_2 - \epsilon_i}(-1) {\bf 1} +
e_{\epsilon_1 + \epsilon_2}(-1)^2 e_{\epsilon_1 - \epsilon_2}(-2)
{\bf 1}, \nonumber \\
&& \quad \mbox{for } i=3, \ldots ,l, \\
&& e_{\epsilon_1 - \epsilon_{2}}(0). e_{\epsilon_1 +
\epsilon_2}(-1)^2 (e_{\epsilon_1 + \epsilon_i}(-1)f_{\epsilon_1 +
\epsilon_i}(-1)+ f_{\epsilon_1 + \epsilon_i}(-1)e_{\epsilon_1 +
\epsilon_i}(-1)){\bf 1} \nonumber \\
&& \ = -2e_{\epsilon_1 + \epsilon_2}(-1)^2 e_{\epsilon_1 +
\epsilon_i}(-1)f_{\epsilon_2 + \epsilon_i}(-1) {\bf 1} +
e_{\epsilon_1 + \epsilon_2}(-1)^2 e_{\epsilon_1 - \epsilon_2}(-2)
{\bf 1}, \nonumber \\
&& \quad \mbox{for } i=3, \ldots ,l, \\
&& e_{\epsilon_1 - \epsilon_{2}}(0). e_{\epsilon_1 +
\epsilon_2}(-1)^2 (e_{\epsilon_1 + \epsilon_2}(-1)f_{\epsilon_1 +
\epsilon_2}(-1)+ f_{\epsilon_1 + \epsilon_2}(-1)e_{\epsilon_1 +
\epsilon_2}(-1)){\bf 1}=0, \\
&& e_{\epsilon_1 - \epsilon_{2}}(0). e_{\epsilon_1 +
\epsilon_2}(-1)^2( e_{\alpha}(-1) f_{\alpha}(-1){\bf 1}
+f_{\alpha}(-1)e_{\alpha}(-1) ){\bf 1}=0, \nonumber \\
&& \quad \mbox{for } \alpha \in \Delta _{D_{l}}^{+} \mbox{ such that
} ( \alpha , \epsilon_1)=0 , ( \alpha , \epsilon_2)=0, \\
&& e_{\epsilon_1 - \epsilon_{2}}(0). e_{\epsilon_1 +
\epsilon_2}(-1)^2 h_{\epsilon_1 - \epsilon_2}(-1)^2 {\bf 1} = -4
e_{\epsilon_1 + \epsilon_2}(-1)^2 e_{\epsilon_1 - \epsilon_2}(-1)
h_{\epsilon_1 - \epsilon_2}(-1) {\bf 1} \nonumber \\
&& \ -4 e_{\epsilon_1 + \epsilon_2}(-1)^2 e_{\epsilon_1 -
\epsilon_2}(-2) {\bf 1}, \\
&& e_{\epsilon_1 - \epsilon_{2}}(0). e_{\epsilon_1 +
\epsilon_2}(-1)^2 h_{\epsilon_1}(-1)^2 {\bf 1} = -4 e_{\epsilon_1 +
\epsilon_2}(-1)^2 e_{\epsilon_1 - \epsilon_2}(-1)
h_{\epsilon_1}(-1) {\bf 1} \nonumber \\
&& \ -4 e_{\epsilon_1 + \epsilon_2}(-1)^2 e_{\epsilon_1 -
\epsilon_2}(-2) {\bf 1}, \\
&& e_{\epsilon_1 - \epsilon_{2}}(0). e_{\epsilon_1 +
\epsilon_2}(-1)^2 h_{\epsilon_1 - \epsilon_2}(-1) h_{\epsilon_1 +
\epsilon_2}(-1){\bf 1} = -2 e_{\epsilon_1 + \epsilon_2}(-1)^2
e_{\epsilon_1 - \epsilon_2}(-1) h_{\epsilon_1 + \epsilon_2}(-1)
{\bf 1}, \nonumber \\
&& \mbox{} \\
&& e_{\epsilon_1 - \epsilon_{2}}(0). e_{\epsilon_1 +
\epsilon_2}(-1)^2 h_{\epsilon_1 - \epsilon_i}(-1) h_{\epsilon_1 +
\epsilon_i}(-1){\bf 1} = - e_{\epsilon_1 + \epsilon_2}(-1)^2
e_{\epsilon_1 - \epsilon_2}(-1) h_{\epsilon_1}(-1) {\bf 1} \nonumber
\\
&& \ -e_{\epsilon_1 + \epsilon_2}(-1)^2 e_{\epsilon_1 -
\epsilon_2}(-2) {\bf 1}, \quad \mbox{for } i=3, \ldots ,l, \\
&& e_{\epsilon_1 - \epsilon_{2}}(0). e_{\epsilon_1 +
\epsilon_2}(-1)^2 h_{\epsilon_1 + \epsilon_2}(-2){\bf 1} =0.
\label{rel.app.zad}
\end{eqnarray}

Using relations (\ref{rel.app.1})-(\ref{rel.app.zad}) and formula
(\ref{rel.sing.v.D}) for vector $v_{D_{l}}$, one directly obtains
$e_{\epsilon_1 - \epsilon_{2}}(0).v_{D_{l}}=0$.

\bibliography{thesis}
\bibliographystyle{plain}

\vskip 1cm

Department of Mathematics, University of Zagreb, Bijeni\v{c}ka 30,
\linebreak 10000 Zagreb, Croatia

E-mail address: perse@math.hr

\end{document}